\documentclass[amssymb,openbib]{article}
\usepackage{amssymb}
\usepackage{epsf}
\usepackage{latexsym,epsf}
\title{A conjecture-generalization of Sondow's formula}
\author{Petros Hadjicostas \\ Department of Mathematics and Statistics \\ Texas Tech University,
Box 41042 \\ Lubbock, TX 79409-1042, USA}
\date{}

\newcommand{\be}{\begin{equation}}
\newcommand{\ee}{\end{equation}}
\newtheorem{lemma}{Lemma}
\newtheorem{theorem}{Theorem}

\newtheorem{corollary}{Corollary}

\newtheorem{conjecture}{Conjecture}
\newcommand{\bl}{\begin{lemma}}
\newcommand{\el}{\end{lemma}}
\newcommand{\bt}{\begin{theorem}}
\newcommand{\et}{\end{theorem}}
\newcommand{\bcor}{\begin{corollary}}
\newcommand{\ecor}{\end{corollary}}
\newcommand{\C}{{\Bbb{C}}}

\begin{document}
\maketitle

\vspace{-.2in}

By generalizing Beukers' arguments in \cite{9} one can easily deduce that (e.g., see \cite{10}): 

\begin{theorem}\label{sinan1}
Let ${\Bbb N}$ be the set of all non-negative integers, and for all $n \in {\Bbb N}\backslash\{0\}$ 
let $d_n$ be the lowest common multiple of $1,2,\ldots,n$.
Let $r, s \in {\Bbb N}$. 

(a) If $r>s$, then for all $n \in {\Bbb N}$,
$$\int_0^1 \int_0^1 \frac{\ln^n(xy)}{1-xy} x^r y^s \; dx\,dy=\frac{n!
(-1)^n}{r-s}\sum_{k=s+1}^r \frac{1}{k^{n+1}},$$
which is a rational number 
whose denominator is a divisor of $d_r^{n+2}$.

(b) If $r=s$, then for all $n \in {\Bbb N}$,
$$\int_0^1 \int_0^1 \frac{\ln^n(xy)}{1-xy} x^r y^r\; dx\, dy=(n+1)!
(-1)^n \left( \zeta(n+2)-\sum_{k=1}^r \frac{1}{k^{n+2}}\right).$$
(Note that the sum of numbers over the empty set is by definition zero.)
\end{theorem}

(For related results see also the papers \cite{14.5}, \cite{9.5}, \cite{8.5}, \cite{11}, \cite{80}, \cite{55}, \cite{giati}, \cite{callie4}, and
\cite{13}. See also the books \cite[pp.\ 365-370]{callie44} and \cite[Appendix A2, pp.\ 372-381]{callie33}.)
Using parts (a) and (b) of Theorem 1, it is easy to show that 
$$\int_0^1 \int_0^1 \frac{\left[-\ln(xy)\right]^{n}}{1-xy} (1-x) \, dx\, dy=\Gamma(n+2)\left[\zeta(n+2)-\frac{1}{n+1}\right],$$
for all $n \in {\Bbb N}$. This gives rise to the following conjecture:

\begin{conjecture}
For $z \in \C$ with $\Re(z) > -2,$
\begin{equation}\label{callie}
\int_0^1 \int_0^1 \frac{\left[-\ln(xy)\right]^{z}}{1-xy} (1-x) \, dx\, dy=\Gamma(z+2)\left[\zeta(z+2)-\frac{1}{z+1}\right].
\end{equation}
\end{conjecture}
The conjecture can {\it probably} be proven by generalizing Beukers' arguments in \cite{9}
and by using Fractional Calculus (e.g., see \cite{333}). 
By taking the limit as $z \rightarrow -1$ on both sides of (\ref{callie}), 
one obtains Sondow's formula
$$\gamma=\int_0^1 \int_0^1 -\frac{1-x}{(1-xy)\ln (xy)}\, dx\, dy,$$
where $\gamma$ is Euler's constant. (See \cite{callie1} and \cite{callie2}.)


\begin{thebibliography}{333}

\bibitem{14.5} Ball, K., and Rivoal, T.\ (2001), ``Irrationalit\'{e} d'une infinit\'{e} de valeurs de la fonction z\^{e}ta aux entiers 
impairs (The irrationality of an infinite number of values of the zeta function at odd integers)," 
{\it Invent.\ Math.}, {\bf 146}, 193--207.

\bibitem{9} Beukers, F.\ (1979), ``A note on the irrationality of $\zeta(2)$
and $\zeta(3)$,'' {\it Bull.\ London Math.\ Soc.}, {\bf 11}, 268-272.

\bibitem{callie44} Borwein, J.M., and Borwein, P.B., {\it Pi and the AGM}, John Wiley, New York, 1987.

\bibitem{callie33} Borwein, P., and Erd\'{e}lyi, T., {\it Polynomials and polynomial inequalities}, Springer, New York, 
1995.

\bibitem{9.5} Dvornicich, R., and Viola, C.\ (1990), ``Some remarks
on Beukers' integrals,"  {\it Number theory}, Vol.~2 (Budapest, 1987),
{\it Colloq.\ Math.\ Soc.\ J\'{a}nos Bolyai}, {\bf 51}, North-Holland,
Amsterdam, 637-657.

\bibitem{10} Hadjicostas, P.\ (2002), ``Some generalizations of Beukers' integrals,'' {\it Kyungpook 
Mathematical Journal}, {\bf 42}, 399-416.

\bibitem{8.5} Hata, M.\ (1995), ``A note on Beukers' integral," {\it J.\ Austral.\ Math.\ Soc.\ Ser.\ A},
{\bf 58}, 143--153.

\bibitem{11} Huylebrouck, D.\ (2001), ``Similarities in irrationality proofs for $\pi$, $\ln 2$, $\zeta(2)$, and
$\zeta(3)$," {\it Amer.\ Math.\ Monthly}, {\bf 108}, 222-231.

\bibitem{80} Krattenthaler, C., and Rivoal, T.\ (2003), ``An identity of Andrews, multiple integrals, and 
very-well-poised hypergeometric series," {\it (http://arXiv.org/abs/math/0312148)},
submitted for publication.


\bibitem{55} Nesterenko, Yu.\ V.\ (2003), ``Integral identities and constructions of approximations to zeta-values," 
Actes des 12\`{e}mes rencontres arithm\'{e}tiques de Caen (June 29--30, 2001), 
{\it J. Th\'{e}orie Nombres Bordeaux}, {\bf 15}, 535--550.

\bibitem{333} Oldham, K.B., and Spanier, J., {\it The Fractional Calculus}, Academic Press, New York, 1974.

\bibitem{giati} Rivoal, T.\  (2000),  ``La fonction z\^{e}ta de Riemann prend
une infinit\'{e} de valeurs irrationnelles aux entiers impairs (There
are infinitely many irrational values of the Riemann zeta function at
odd integers)," {\it C.\ R.\ Acad.\ Sci.\ Paris S\'{e}r.\ I Math.}, {\bf 331}, 267-270.

\bibitem{callie1} Sondow, J.\ (2004), ``A faster product for $\pi$ and a new integral for $\ln(\pi/2)$," 
to appear in {\it The American Mathematical Monthly}, {\it (http://arXiv.org/abs/math.NT/0401406).}

\bibitem{callie2} Sondow, J.\ (2004), ``Double integrals for Euler's constant and $\ln(4/\pi)$," submitted to the {\it The 
American Mathematical Monthly}, {\it (http://arXiv.org/abs/math.CA/0211148).}

\bibitem{callie4} Zudilin, W.\ (2002), ``Very well-poised hypergeometric series and multiple integrals," 
{\it Russian Math.\ Surveys}, {\bf 57}, 824--826.

\bibitem{13} Zudilin, W.\ (2004), ``Well-poised hypergeometric transformations of Euler-type multiple integrals," 
{\it J.\ London Math.\ Soc.\ (2)}, {\bf 70}, 16 pages. 


	
\end{thebibliography}
\end{document}